\documentclass{amsart}
\usepackage{epsf}
\def\hepsffile{\leavevmode\epsffile}

\theoremstyle{plain}
\newtheorem{thm}{Theorem}[subsection]
\newtheorem{thm1}{Theorem}

\newtheorem{lem}[thm]{Lemma}

\newtheorem{prop}[thm]{Proposition}

\theoremstyle{definition}
\newtheorem{defin}[thm]{Definition}

\newtheorem{emf}[thm]{}

\def\id{\protect\operatorname{id}}

\def\sign{\protect\operatorname{sign}}

\def\pr{\protect\operatorname{pr}}

\def\C{{\mathbb C}}
\def\Z{{\mathbb Z}}

\def\R{{\mathbb R}}
\def\N{{\mathbb N}}

\def\1{\hbox{\rm\rlap {1}\hskip.03in{\rom I}}}
\def\Bbbone{{\rm1\mathchoice{\kern-0.25em}{\kern-0.25em}
	{\kern-0.2em}{\kern-0.2em}I}}

\def\p{\partial}

\begin{document}
\hyphenation{Ca-m-po}
\title[Vassiliev Invariants of Legendrian and Pseudo Legendrian Knots]
{Isomorphism of the groups of 
Vassiliev Invariants of Legendrian and of Pseudo Legendrian Knots 
in Contact $3$-manifolds}
\author[V.~Tchernov (Chernov)]{Vladimir Tchernov (Chernov)}
\address{Institute for Mathematics, Zurich University, Winterthurerstrasse
190, CH-8057, Zurich, Switzerland}
\email{Chernov@math.unizh.ch}
\begin{abstract}
The study of the Vassiliev invariants of Legendrian knots was started by
D.~Fuchs and S.~Tabachnikov who showed that the groups of $\C$-valued
Vassiliev invariants of Legendrian and of framed knots in the standard
contact $\R^3$ are canonically isomorphic.

Recently we constructed the first examples where Vassiliev
invariants of Legendrian and of framed knots are different, and Vassiliev
invariants of Legendrian knots distinguish
Legendrian knots that are isotopic as framed knots and homotopic as
Legendrian immersions. This raised the question what information about
Legendrian knots can be captured using Vassiliev invariants.

Here we answer this question by showing that for any contact $3$-manifold
with a cooriented contact structure the groups of Vassiliev invariants of
Legendrian knots and of knots that are nowhere tangent to a vector field
that coorients the contact structure are canonically isomorphic.
\end{abstract}
\maketitle

\section{Introduction} 
In this section we describe the main results of the paper.
(In case any of the terminology appears to be new to the reader, the
corresponding definitions are given in the next section.)

If a contact structure on a $3$-manifold is cooriented, then every
Legendrian knot (i.e. a knot that is everywhere tangent to the contact
distribution) has a natural framing (a continuous normal vector field).
Hence when studying Legendrian knots in such contact manifolds the main
question is to distinguish those of them that realize isotopic framed knots.

On the other hand a cooriented contact structure $C$ on a manifold $M$ gives 
rise to a nondegenerate vector field $V_C$ in $TM$ that coorients the
contact structure.
Clearly if two Legendrian knots $K_0$ and $K_1$ are
isotopic as Legendrian knots, then they are also isotopic in the category of 
knots that are everywhere nontangent to $V_C$.

This observation leads to the following definition. 
Let $(M,C)$ be a contact manifold with a cooriented contact structure, 
and let $V_C$ be the nondegenerate vector field on $TM$ that coorients the
contact structure. A knot $K_0$ in $M$ is said to be {\em
Pseudo-Legendrian\/} if it is everywhere nontangent to $V_C$.
Two Pseudo-Legendrian knots $K_0$ and $K_1$ (in $(M,C,V_C)$)
are {\em  Pseudo-Legendrian isotopic\/} if there
exists an isotopy $I:[0,1]\times S^1\rightarrow M$ such that 
$I\big|_{0\times S^1}=K_0$, $I\big|_{1\times S^1}=K_1$, and $\forall t\in
[0,1]$ the knot $K_t=I\big|_{t\times S^1}$ is Pseudo-Legendrian (with
respect to $V_C$).

Clearly if $K_0$ and $K_1$ are Legendrian isotopic Legendrian knots, 
then they are also Pseudo-Legendrian isotopic (with respect to any
nondegenerate vector field $V_C$ that coorients $C$). 

Vassiliev invariants proved to be an extremely useful tool in the study of
framed knots, and the conjecture is that they are sufficient to distinguish
all the isotopy classes of framed knots. Vassiliev invariants can also be 
easily defined in the categories of Legendrian and of Pseudo-Legendrian knots.

The study of the groups of Vassiliev invariants of Legendrian knots was
initiated by the work~\cite{FuchsTabachnikov} of D.~Fuchs and S.~Tabachnikov
where it was proved that the groups of $\C$-valued Vassiliev invariants of 
Legendrian and of framed knots in the standard contact $\R^3$ are
canonically isomorphic. Later the similar result was proved by
J.~Hill~\cite{Hill} for
the groups of $\C$-valued Vassiliev invariants of Legendrian and of framed
knots in the spherical cotangent bundle $ST^*\R^2$ of $\R^2$ with the
standard contact structure. The proofs of these isomorphisms were based on
the existence of the universal $\C$-valued Vassiliev invariant for these
spaces, also known as the Kontsevich integral~\cite{Kontsevich}. (For
$ST^*\R^2$ such universal invariant was first constructed by
V.~Goryunov~\cite{Goryunov}.) Unfortunately the Kontsevich integral exists
only for a rather limited collection of $3$-manifolds. (Recently Andersen,
Mattes, and Reshetikhin~\cite{AMR} constructed such an invariant for manifolds that are
$\R^1$-fibered over an oriented surface $F$ with $\p F\neq \emptyset$.) For
this reason the question whether the groups of Vassiliev invariants of
Legendrian and of framed knots are always isomorphic was open for some time.

Recently the author used different technique to
prove~\cite{Chernovpreprint},~\cite{VasLegendrian} 
that for any Abelian group $\mathcal
A$ the groups of $\mathcal A$-valued Vassiliev invariants of Legendrian and of 
framed knots are canonically isomorphic for a large class of 
contact $3$-manifolds with a cooriented contact structure. This class of
contact $3$-manifolds $(M,C)$ includes all contact manifolds with a tight contact
structure, all contact manifolds that are closed and admit a metric of
negative sectional curvature, and all contact manifolds such that the Euler
class of the contact bundle is in the torsion of $H^2(M,\Z)$.

On the other
hand~\cite{Chernovpreprint},~\cite{VasLegendrian}
the author constructed the first known examples of contact manifolds
where the groups of Vassiliev invariants of Legendrian 
and of framed knots are not
canonically isomorphic. In these examples Vassiliev invariants of Legendrian
knots can be successfully used to distinguish Legendrian knots that realize
isotopic framed knots and that are homotopic  as Legendrian immersions of $S^1$.
Namely, such examples were constructed for $M=S^1\times S^2$ and for any $M$ that
is an orientable total space of an $S^1$-bundle over a nonorientable surface
of genus bigger than one. This brought up a question what information about
Legendrian knots can be captured with the help of Vassiliev invariants of
Legendrian knots.

\def\mfootu{\footnote{
In their work~\cite{BenedettiPetronio}, p. 34, R.~Benedetti and C.~Petronio
conjectured the fact that is very similar to the one shown in
Theorem~\ref{first}, but their definition of
Pseudo-Legendrian isotopy is different. Namely, let $\mathcal V$ be the space
of nondegenerate vector fields on $TM$ that are homotopic (as nondegenerate
vector fields) to a vector field that coorients $C$.
They call a Pseudo-Legendrian knot
in $(M,C)$ a pair $(K, V)$ that consists of $V\in\mathcal V$ and of a knot
$K$ that is everywhere nontangent to $V$. 
They say that Pseudo-Legendrian knots $(K_0, V_0)$ and $(K_1, V_1)$ 
are Pseudo-Legendrian isotopic if there exists a homotopy
$I_V:[0,1]\rightarrow \mathcal V$ with $I_V(0)=V_0$, $I_V(1)=V_1$, and an
isotopy $I:[0,1]\times S^1\rightarrow M$ with $I\big|_{0\times
S^1}=K_0$, $I\big|_{1\times S^1}=K_1$ such that $\forall t\in [0,1]$ 
$K_t=I\big|_{t\times S^1}$ is nowhere tangent to $I_V(t)\in\mathcal V$.
We were not able to prove their conjecture and as it seems to us it should 
be possible to construct an example showing that the groups of Vassiliev
invariants of Legendrian knots and of knots 
that are Pseudo-Legendrian with respect to their
definition are different.}}

Here we answer this question\protect\mfootu by proving the following Theorem, that says
that the groups of Vassiliev invariants of Legendrian and of
Pseudo-Legendrian knots are always canonically isomorphic.

\let\mfootu\relax

Let $\mathcal A$ be an Abelian group, 
let $(M,C)$ be a contact $3$-manifold with a cooriented contact structure,
and let $V_C$ be a nondegenerate vector field that coorients $C$. 
Let $\mathcal L$ be a connected component of the space of Legendrian
immersions of $S^1$ and let $\mathcal L_p$ be a connected component of the
space of Pseudo-Legendrian immersions of $S^1$ (with respect to $V_C$) that
contains $\mathcal L$.
(Such a component always exists since a path in $\mathcal L$ corresponds to
a path in $\mathcal L_p$. Moreover one can show, see
Proposition~\ref{components}, that every component of the space of
Pseudo-Legendrian immersions of $S^1$ contains a unique component of the
space of Legendrian immersions of $S^1$.)
\begin{thm1}\label{first}
The groups of $\mathcal A$-valued
Vassiliev invariants of Legendrian knots from $\mathcal L$ and of
Pseudo-Legendrian knots from $\mathcal L_p$ are canonically isomorphic.
\end{thm1}
See Theorem~\ref{isomorphism}.

In particular, if
$K_1$ and $K_2$ are two Legendrian knots that are homotopic as Legendrian
immersions and that realize isotopic framed knots, and $x$ is a Vassiliev
invariant of Legendrian knots such that $x(K_1)\neq x(K_2)$, then $K_1$ and
$K_2$ are not isotopic as Pseudo-Legendrian knots. This means that the only
information about a Legendrian knot that can be
captured using Vassiliev invariants of Legendrian knots is the
Pseudo-Legendrian isotopy class of the Legendrian knot.

\section{Conventions and definitions.}\label{definitions}
We work in the smooth category.

In this paper $\mathcal A$ is an Abelian group
(not necessarily torsion free), and $M$ is a
connected oriented $3$-dimensional 
Riemannian 
manifold (not necessarily
compact).

A {\em contact structure\/} on a $3$-dimensional manifold $M$ is a smooth
field $\{ C_x\subset T_xM|x\in M\}$ of tangent $2$-dimensional planes,
locally defined as a kernel of a differential $1$-form $\alpha$ with
non-vanishing
$\alpha\wedge d\alpha$. A manifold with a contact structure possesses
the canonical orientation determined by the volume form $\alpha\wedge\ d
\alpha$. The standard contact structure in $\R^3$ is the kernel of the
$1$-form $\alpha=ydx-dz$.

A {\em contact element\/} on a manifold is a hyperplane
in the tangent space to the manifold at a point.
For a surface $F$ we denote by $ST^*F$ the space of all
cooriented (transversally oriented) contact elements of $F$. This space is
the spherical cotangent bundle of $F$ and it is one of the classical
examples of contact manifolds. Its natural contact structure is the
distribution of tangent hyperplanes given by the condition
that the velocity vector of the incidence point of a contact element
belongs to the element.

A contact structure is {\em cooriented\/} if the
$2$-dimensional planes defining the contact structure are continuously
cooriented (transversally oriented).
A contact structure is {\em oriented\/} if the $2$-dimensional planes
defining the
contact structure are continuously oriented. Since every contact manifold
has a natural orientation we see that every cooriented contact structure is
naturally oriented and every oriented contact structure is naturally
cooriented.

A contact structure is {\em
parallelizable\/} ({\em parallelized\/}) if the $2$-dimensional vector
bundle
$\{ C_x \}$ over $M$ is trivializable (trivialized). Since every contact
manifold
has a canonical orientation, one can see that every parallelized contact
structure is naturally cooriented. 

A contact structure $C$ on a manifold $M$
is said to be {\em overtwisted\/} if there exists a $2$-disk $D$ embedded
into $M$ such that the boundary $\p D$ is tangent to $C$ while the disk $D$
is transverse to $C$ along $\p D$. Not overtwisted contact structures are
called {\em tight\/}.

A {\em curve\/} in $M$ is an immersion of $S^1$ into $M$.
(All curves have the natural orientation induced by the orientation of
$S^1$.) A {\em framed curve\/} in $M$ is a curve equipped
with a continuous unit normal vector field.
 
A {\em Legendrian curve\/} in a contact manifold $(M,C)$ is a curve
in $M$ that is everywhere tangent to $C$. If the contact structure on
$M$ is cooriented, then every Legendrian curve has a natural framing given
by the unit normals to the planes of the contact structure that point in the
direction specified by the coorientation.
 
To a Legendrian curve $K_l$ in a contact manifold $(M,C)$ 
with a parallelized contact structure one can associate an integer that is the
number of revolutions of the direction of the velocity vector of $K_l$ (with
respect to the chosen frames in $C$) under traversing $K_l$ according
to the orientation. This integer is called the {\em Maslov number\/} of
$K_l$. The set of Maslov numbers enumerates the
set of the connected components of the space of Legendrian
curves in $\R^3$ (cf.~\ref{h-principleLegendrian}).

For a contact manifold $(M,C)$ with a cooriented contact structure fix 
a nondegenerate vector field $V_C$ that coorients the contact structure.
A {\em Pseudo-Legendrian curve\/} in $(M,C, V_C)$ is a curve that is nowhere
tangent to $V_C$.
Clearly every Legendrian curve in $(M,C)$ {\em realizes a Pseudo-Legendrian
curve.\/} (This means that if $L$ is a Legendrian curve in $(M,C)$, 
then it is also a Pseudo-Legendrian curve in $(M,C,V_C)$.)

A {\em knot (framed knot)\/} in $M$ is an embedding (framed embedding) of
$S^1$ into $M$.
In a similar way we define Legendrian knots, and Pseudo-Legendrian knots in
a contact manifold $(M,C)$ with a cooriented contact structure.

A {\em singular (framed)\/} knot with $n$ double points is a curve (framed
curve)
in $M$ whose only singularities are $n$ transverse double points.
An {\em isotopy\/} of a singular (framed) knot
with $n$ double points is a path in the space of singular (framed) knots
with $n$ double points under which the preimages of the double points on $S^1$
change continuously. In a similar way we define singular Legendrian and
Pseudo-Legendrian knots and the notion of isotopy of singular Legendrian
knots and of singular Pseudo-Legendrian knots.
 
An $\mathcal A$-valued framed (resp. Legendrian, resp. Pseudo-Legendrian) knot
invariant is an $\mathcal A$-valued function on the set of the isotopy classes 
of framed (resp. Legendrian, resp. Pseudo-Legendrian) knots.

A transverse double point $t$ of a singular knot can be resolved in two
essentially different ways. We say that a resolution of a double point is
positive (resp. negative) if the tangent vector to the
first strand, the tangent vector to the second strand, and the vector from
the second strand to the first form the positive $3$-frame. (This does
not depend on the order of the strands).
If the singular knot is Legendrian (resp. Pseudo-Legendrian),
then these resolution can be made in the
category of Legendrian (resp. Pseudo-Legendrian) knots.
 
A singular framed (resp. Legendrian, resp. Pseudo-Legendrian) knot $K$ with $(n+1)$
transverse double points
admits $2^{n+1}$ possible resolutions of the double points. The sign of the
resolution
is put to be $+$ if the number of negatively resolved double points is even,
and
it is put to be $-$ otherwise.
Let $x$ be an $\mathcal A$-valued invariant of framed (resp. Legendrian,
resp. Pseudo-Legendrian) knots. The invariant $x$ is said to be of {\em finite
order\/} (or {\em Vassiliev invariant\/}) if there exists a nonnegative
integer $n$ such that for any singular knot $K_s$ with $(n+1)$
transverse double points the sum (with appropriate signs) of the values of
$x$ on the nonsingular
knots obtained by the $2^{n+1}$ resolutions of the double points is zero.
An invariant is said to be of order not greater than $n$ (of order $\leq n$)
if $n$
can be chosen as the integer in the definition above. The group of $\mathcal
A$-valued finite order invariants has an increasing filtration by the
subgroups of the invariants of order $\leq n$.
 
\begin{emf}\label{h-principleLegendrian}
{\em $h$-principle for Legendrian curves.\/}
For $(M,C)$ a contact manifold with a cooriented contact structure, we put
$CM$ to be the total space of the fiberwise spherization of the contact
bundle, and we put $\pr:CM\rightarrow M$ to be the corresponding locally
trivial $S^1$-fibration. The $h$-principle proved for the Legendrian curves
by M.~Gromov (\cite{Gromov}, pp.338-339) says that the space of Legendrian
curves in
$(M,C)$ is weak homotopy equivalent to the space of free loops
$\Omega CM$
in $CM$. The equivalence is given by mapping a point of a Legendrian curve
to the point of $CM$ corresponding to the direction of the velocity vector
of the curve at this point. In particular the $h$-principle implies that the
set of the connected components of the space of Legendrian curves can be
naturally identified with the set of the conjugacy classes of elements of
$\pi_1(CM)$.
\end{emf}
 
\begin{emf}\label{description}{\em Description of Legendrian and
of Pseudo-Legendrian knots in $\R^3$.\/}
The contact Darboux theorem says that every contact $3$-manifold $(M,C)$ is
locally contactomorphic to
$\R^3$ with the standard contact structure that is the kernel of the
$1$-form
$\alpha=ydx-dz$. A chart in which $(M,C)$ is contactomorphic to the standard
contact $\R^3$ is called {\em a Darboux chart.\/}
 
Legendrian knots in the standard contact $\R^3$
are conveniently presented by the projections
into the plane $(x,z)$. Identify a  point $(x, y, z)\in \R^3$ with the
point $(x,z)\in \R^2$ furnished with the fixed
direction of an unoriented
straight line through $(x,z)$ with the slope $y$. Then the curve in $\R^3$
is a one parameter family of points with non-vertical directions in $\R^2$.

While a generic regular curve has a regular projection into the
$(x,z)$-plane, the projection of a generic Legendrian curve into the
$(x,z)$-plane has isolated critical points (since all the planes of the
contact structure are parallel to the $y$-axis). Hence the projection of a
generic Legendrian curve may have cusps. A curve in $\R^3$ is Legendrian if
and only if the corresponding planar curve with cusps
is everywhere tangent to the field of directions. In particular this field
is determined by the curve with cusps.

A Pseudo-Legendrian knot in $(\R^3, \ker (ydx-dz))$ can be depicted 
as follows. 
Let $V_C$ be a unit vector field on $\R^3$ that coorients the contact
structure.
Choose a system of coordinates $(x',y',z')$ in $\R^3$ 
so that $V_C$ is parallel to the $z'$-axis and points in the
same direction. 
Then a Pseudo-Legendrian knot $K$ can be depicted 
by the standard knot diagram in $(x',y',z')$. 
Since $K$ is Pseudo-Legendrian it means that the velocity
vector of $K$ at every point is not parallel to the $z'$-axis. A
Pseudo-Legendrian isotopy of a Pseudo-Legendrian knot
can be depicted by a sequence of second and third
Reidemeister moves. (The first Reidemeister move does not occur, since
during this move the
velocity vector of one of the points on the kink becomes parallel to the
$z'$-axis.) 
\end{emf}

\subsection{Isomorphism between the groups of order $\leq n$ invariants of
Legendrian and of Pseudo-Legendrian knots}

Let $(M,C)$ be a contact manifold with a cooriented contact structure,
let $V_C$ be a nondegenerate vector field that coorients the contact
structure, 
let $\mathcal L$ be a connected component of
the space of Legendrian curves in $(M,C)$, and let $\mathcal L_p$ be the
connected component of the space of Pseudo-Legendrian curves 
in $(M,C,V_C)$ that contains $\mathcal L$.
(Such a component exists because a Legendrian curve $L$ in $(M,C)$ is 
Pseudo-Legendrian in $(M,C,V_C)$. Moreover as it is shown in
Proposition~\ref{components} every component of the space of
Pseudo-Legendrian curves in $(M,C,V_C)$ contains a unique component of the
space of Legendrian curves in $(M,C)$.)
Let $V_n^{\mathcal L}$ (resp. $V_n^{\mathcal L_p}$)
be the group of $\mathcal A$-valued order $\leq n$ invariants of Legendrian
(resp. Pseudo-Legendrian) knots from $\mathcal L$ (resp.
from $\mathcal L_p$). Clearly every invariant $y\in V_n^{\mathcal L_p}$
restricted to
the category of Legendrian knots in $\mathcal L$ is an element $\phi(y)\in
V_n^{\mathcal L}$. This gives a homomorphism $\phi:V_n^{\mathcal
L_p}\rightarrow V_n^{\mathcal L}$.

\begin{thm}\label{donotdistinguish}
$x(K_1)=x(K_2)$,
for every $x\in V_n^{\mathcal L}$ and for every Legendrian knots
$K_1, K_2\in\mathcal L$ such that $K_1$ and $K_2$ are
Pseudo-Legendrian isotopic knots in $(M,C,V_C)$.
\end{thm}

For the Proof of Theorem~\ref{donotdistinguish} see
Section~\ref{proofdonotdistinguish}.

\begin{thm}\label{equivalent}
The following two statements {\textrm I} and \textrm{II} are equivalent.
\begin{description}
\item[\textrm{I}] $\phi:V_n^{\mathcal L_p}\rightarrow V_n^{\mathcal L}$ is an
isomorphism.
\item[\textrm{II}] $x(K_1)=x(K_2)$ for every
$x\in V_n^{\mathcal L}$
and for every Legendrian knots
$K_1, K_2\in\mathcal L$ such that $K_1$ and $K_2$
are Pseudo-Legendrian isotopic knots in $(M,C,V_C)$.
\end{description} 
\end{thm}

The Proof of Theorem~\ref{equivalent} becomes obvious when
the mapping from the Legendrian isotopy classes of Legendrian knot from
$\mathcal L$ to the Pseudo-Legendrian isotopy classes of Pseudo-Legendrian
knots from $\mathcal L_p$ is surjective. However the famous Bennequin inequality 
shows that this map is not surjective even when $(M,C)$ is the standard
contact $\R^3$.

The Proof of Theorem~\ref{equivalent} is given in
Section~\ref{proofequivalent}.

Combining Theorems~\ref{donotdistinguish} and~\ref{equivalent} we get 
the following.
\begin{thm}\label{isomorphism}
The groups $V_n^{\mathcal L}$ and $V_n^{\mathcal L_p}$ 
of $\mathcal A$-valued Vassiliev invariants of Legendrian knots from
$\mathcal L$ and of Pseudo-Legendrian knots from $\mathcal L_p$ are
canonically isomorphic. 
\end{thm}

\section{Proofs}
\subsection{Proof of
Theorem~\ref{donotdistinguish}}\label{proofdonotdistinguish}
\begin{emf}{\em Some Useful facts proved by D.~Fuchs and
S.~Tabachnikov~\cite{FuchsTabachnikov}.}
There are two types of cusps arising under the projection of the
part of a Legendrian knot that is contained in a Darboux chart to the
$(x,z)$-plane
(see~\ref{description}). They are formed by cusps for which the
branch of the projection of the knot going away from the cusp is locally
located respectively above or below the tangent line at the cusp point. (See
Figure~\ref{twocusp.fig}~a and~b respectively.) For a Legendrian knot $K$ and $i,j\in\N$ we denote
by
$K^{-i,-j}$ the Legendrian knot
obtained from $K$ by the modification corresponding to an addition of
$i$ cusp pairs of the first type and $j$ cusp pairs of the second type to
the projection of the part of $K$ located in a Darboux chart.

\begin{figure}[htbp]
 \begin{center}
  \epsfxsize 8cm
  \hepsffile{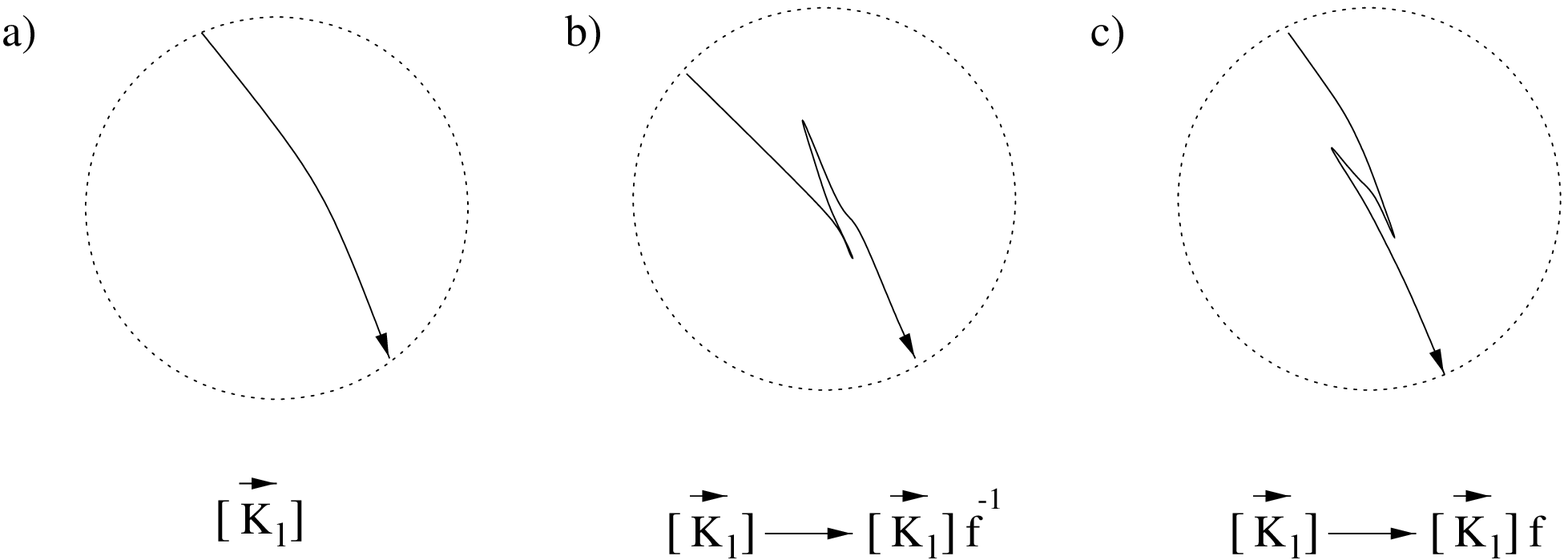}
 \end{center}
\caption{}\label{twocusp.fig}
\end{figure}

The following three facts were proved by Fuchs and Tabachnikov~\cite{FuchsTabachnikov}. 	

\begin{description}
\item[1] Let $K_1$ and $K_2$ be Legendrian knots in
the standard contact $\R^3$ that realize isotopic unframed knots. Then for
any $n_1, n_2\in\N$ large enough there exist $n_3, n_4\in \N$ such that
the Legendrian knot $K_1^{-n_1,-n_2}$ is Legendrian isotopic to
$K_2^{-n_3,-n_4}$.
\item[2] If there exists $n\in\N$ such that Legendrian knots $K_1^{-n,-n}$ and
$K_2^{-n,-n}$ are Legendrian isotopic, then every Vassiliev invariant of
Legendrian
knots takes equal values on $K_1$ and on $K_2$.
\item[3] The number $n$ from the previous observation exists if the ambient
contact manifold is $\R^3$ and the Legendrian knots $K_1$ and $K_2$ belong
to the
same component of the space of Legendrian curves and realize isotopic framed
knots.
\end{description}

The first two observation are true for any contact $3$-manifold (since the
proof
of the corresponding facts is local). But the number $n$ from the
statement
of the third observation does not exist in general.
In the case of the
ambient manifold being $\R^3$ Fuchs and Tabachnikov showed the existence of
such $n$ using the explicit calculation involving the Maslov classes and
Bennequin invariants of Legendrian knots. However in order for the
Bennequin invariant to be well-defined the
knots have to be zero-homologous, and in order for the Maslov class to be
well-defined the knots have to be zero-homologous or the contact structure
has to be parallelizable.
\end{emf}

Clearly to prove Theorem~\ref{donotdistinguish} it suffices to show that the number
$n$ from the third observation exists for any
Legendrian knots $K_1$ and $K_2$ that realize isotopic Pseudo-Legendrian knots
in $(M,C,V_C)$.
(We assume that the contact structure $C$ on $M$ is cooriented.)

Let $K_1$ and $K_2$ be Legendrian knots that realize Pseudo-Legendrian
isotopic Pseudo-Legendrian knots (in $(M,C,V_C)$), and
let $n_1, n_2, n_3, n_4\in \N$ be such that $K_1^{-n_1, -n_2}$ and
$K_2^{-n_3,-n_4}$ are Legendrian isotopic.

We start by showing that if $K_1$ and $K_2$ are Pseudo-Legendrian isotopic
then $n_1, n_2, n_3, n_4$ can be chosen so that
$n_1+n_2=n_3+n_4$, and that
$n_1-n_2=n_3-n_4$.

\begin{emf}{\em Proof of the fact that $n_1, n_2, n_3, n_4$ can be chosen so
that $n_1+n_2=n_3+n_4$.}
Let $I:S^1\times [0,1]\rightarrow M$ be
the Pseudo-Legendrian isotopy that changes $K_1$ to $K_2$.

Analyzing the proof of Fuchs and Tabachnikov one verifies that for $n_1, n_2$ 
large enough the Legendrian isotopy $\mu:S^1\times[0,1]\rightarrow M$
changing $K_1^{-n_1, -n_2}$ to $K_2^{-n_3, -n_4}$ can be
chosen so that for every $t\in [0,1]$ the Legendrian knot 
$\mu^t=\mu\big|_{(S^1\times t)}$ is contained in a thin tubular neighborhood
$T_t$ of $I_t=I\big|_{(S^1\times t)}$ and is isotopic (as an
unframed knot) to $\mu ^t$ inside $T_t$.

Every Pseudo-Legendrian knot is naturally framed. Put $\bar \mu^t$ and $\bar
I_t$ to be the framed knots corresponding respectively to $\mu^t$ and to $I_t$
For two framed knots $\bar \mu^t$ and  
$\bar I_t$ realizing unframed knots
that are isotopic inside $T_t$ there is a well-defined $\Z$-valued
obstruction to be isotopic inside $T_t$ in the category of framed knots.

This obstruction is the difference of the self-linking numbers of the
inclusions of $\bar \mu ^t$ and $\bar I_t$ into $\R^3$ induced by an
identification of $T_t$ with the standard solid torus in $\R^3$. (One
verifies that for  $\bar \mu ^t$ and $\bar I_t$ that are isotopic as unframed
knots inside $T_t$ this difference does not depend on the choice of the
identification of $T_t$ with the standard solid torus in $\R^3$.)

From the formula for the Bennequin invariant stated
in~\cite{FuchsTabachnikov}
one gets that the value of the obstruction for $K_1^{-n_1, -n_2}$ to be
isotopic as a framed knot to $K_1$ inside $T_0$ is equal to $n_1+n_2$.
Similarly the value of the obstruction for $K_2^{-n_3, -n_4}$ to be
isotopic as a framed knot to $K_2$ inside $T_1$ is equal to $n_3+n_4$.
Clearly the value of the obstruction for $\bar \mu^t$ to be isotopic to
$I_t$ inside $T_t$ does not depend on $t$, and  we get that $n_1+n_2=n_3+n_4$.
\end{emf}

\begin{emf}{\em Proof of the fact that $n_1, n_2, n_3, n_4$ can be chosen so
that $n_1-n_2=n_3-n_4$.}

Identify $S^1$ with the circle in $\C$ that consists of the numbers of the
absolute value one. Put $[1,i]$, $[i,-1]$, $[-1, -i]$, and $[-i, 1]$ to be
the four nonoverlapping arcs of $S^1$ with the end points at $1, -1, i, -i\in
S^1$.

Let $p:S^1\times S^1\rightarrow M$ be the mapping, such that
\begin{description}
\item[1] $p\big|_{S^1\times [1,i]}$ is a homotopy of loops that
changes $K_1$ to $K_1^{-n_1, -n_2}$ and happens in a thin tubular
neighborhood of $K_1$.
\item[2] $p\big|_{S^1\times [i,-1]}$ is (up to the reparametrization) the
Legendrian isotopy $\mu$ that changes $K_1^{-n_1, -n_2}$ to $K_2^{-n_3, -n_4}$.
\item[3] $p\big|_{S^1\times [-1, -i]}$ is a homotopy of loops that changes
$K_2^{-n_3, -n_4}$ to $K_2$ and happens in a thin tubular neighborhood of
$K_2$.
\item[4] $p\big|_{S^1\times [-i, 1]}$ is (up to the reparametrization) the
isotopy $I^{-1}$ that changes $K_2$ to $K_1$.
\end{description}

Consider the oriented $\R^2$ bundle $\pr:\xi\rightarrow S^1\times S^1$ that is induced
by $p:S^1\times S^1\rightarrow M$ from the contact bundle $C$ on $M$. Since the
isotopies $I_K$ and $\mu$ were chosen so that they are $C^0$
close, we get that the homology class realized by 
$p:S^1\times S^1\rightarrow M$ 
in $H_2(M)$ is zero. Thus the Euler class $e_\xi$ of $\pr:\xi\rightarrow
S^1\times S^1$ is $0\in\Z=H^2(S^1\times S^1)$.

Put $\pr_1:\xi_1\rightarrow S^1\times [1,i]$, $\pr_2:\xi_2\rightarrow
S^1\times [i,-1]$, $\pr_3:\xi_3\rightarrow S^1\times [-1,-i]$, and 
$\pr_4:\xi_4\rightarrow S^1\times [-i, 1]$ to be the restrictions of
$\pr:\xi\rightarrow S^1\times S^1$.

The velocity vectors of Legendrian knots 
$K_1$, $K_1^{-n_1, -n_2}$, $K_2^{-n_3, -n_4}$, and $K_2$ give the section of
$\pr:\xi\rightarrow S^1\times S^1$ over $\p (S^1\times [1,i])$, 
$\p (S^1\times [i,-1])$, $\p(S^1\times [-1, -i]),$ and 
$\p(S^1\times [-i, 1])$. 

Put $e_{\xi_1}$ to be the $\Z$-valued obstruction to extend the nonzero section of
$\pr_1:\xi_1\rightarrow S^1\times [1,i]$ over 
$\p (S^1\times [1,i])=S^1\times \{1\}\cup S^1\times \{i\}$ given by the
velocity vectors of $K_1$ and of $K_1^{-n_1, -n_2}$ to the nonzero section of
$\pr_1:\xi_1\rightarrow S^1\times [1,i]$ over $S^1\times [1,i]$. (It is the
relative Euler class that takes values in $\Z=H^2(S^1\times [1,i],
\p (S^1\times [1,i]))$.)

Similarly put $e_{\xi_2},e_{\xi_3}, e_{\xi_4}\in\Z$ to be the obstructions
to extend the nonzero sections of $\xi_2, \xi_3, \xi_4$ over 
$\p (S^1\times
[i,-1]), \p (S^1\times [-1, -i]),$ and $\p (S^1\times [-i, 1])$ that were described
above to the nonzero sections respectively over $S^1\times                                       
[i,-1], S^1\times [-1, -i],$ and $S^1\times [-i, 1]$.

One verifies that 
\begin{equation}\label{sume}
e_{\xi_1}+e_{\xi_2}+e_{\xi_3}+e_{\xi_4}=e_{\xi}=0.
\end{equation} 

For a Legendrian knot $K$ in the standard contact $\R^3$ put
$m(K)\in\Z$ to be the Maslov class of $K$.
In~\cite{FuchsTabachnikov} it is shown that $m(K^{-i_1,
-i_2})=m(K)+(i_2-i_1)$ for $i_1,i_2\in\Z$ and for the
Legendrian knots $K$ and $K^{-i_1, -i_2}$ in the standard contact 
$\R^3$.
A straightforward verification based on this equality 
shows that 
\begin{equation}\label{exi1exi3}
e_{\xi_1}=n_2-n_1,\text{ and }e_{\xi_3}=n_3-n_4. 
\end{equation}

The velocity vectors of the Legendrian knots $\mu(t)$ define the nonzero 
section of
$\pr_2:\xi_2\rightarrow S^1\times [i,-1]$ that extends the nonzero 
section over $\p
(S^1\times [i,-1])$ defined by the velocity vectors of $K_1^{-n_1, -n_2}$ and of
$K_2^{-n_3, -n_4}$. Thus 
\begin{equation}\label{exi2}
e_{\xi_2}=0.
\end{equation}

Since $\forall t\in [0,1]$ the knot $I_t=I\big|_{S^1\times t}$ is nowhere tangent to $V_C$ we get
that for every point of $I_t$ the projection to $C$ along $V_C$ of the velocity
vector of $I_t$ at this point is nonzero. Thus the isotopy $I$ defines the
nonzero section of $\xi_4$ over $S^1\times [-i,1]$ that is the extension 
of the nonzero section of $\xi_4$ over $\p (S^1\times [-i,1])$ defined 
by the velocity vectors of $K_1$ and of $K_2$. Hence we get that 
\begin{equation}\label{exi4}
e_{\xi_4}=0.
\end{equation}

Combining together identities~\eqref{sume},~\eqref{exi1exi3},~\eqref{exi2},
and~\eqref{exi4}, we get that $0=e_{\xi_1}+e_{\xi_2}+e_{\xi_3}+e_{\xi_4}=
e_{\xi_1}+e_{\xi_3}=(n_2-n_1)+(n_3-n_4)$. Thus $n_1-n_2=n_3-n_4$.
\end{emf}

\begin{emf}
From the identities $n_1+n_2=n_3+n_4$ and $n_1-n_2=n_3-n_4$ one gets that
$n_1=n_3$ and $n_2=n_4$. Assume that $n_1\geq n_2$. (The case where
$n_2>n_1$ is treated similarly.) Put $k=n_1-n_2$.
It is easy to show that
since $K_1^{-n_1, -n_2}$ and $K_2^{-n_3,-n_4}$ are Legendrian isotopic,
then $K_1^{-n_1, -n_2-k}$ and $K_2^{-n_3,-n_4-k}$ are also Legendrian isotopic.
(Basically one can keep the $k$ extra cusp pairs
close together on a small piece of the projection of the part of
the knot contained in a Darboux chart during the whole isotopy
process.) But $K_1^{-n_1, -n_2-k}$ and $K_2^{-n_3,-n_4-k}$ are obtained from
$K_1$ and $K_2$
by the
modification corresponding to the addition of $n_1=n_2+k=n_3=n_4+k$ pairs of
cusps of each of the two types, and we can take $n$ from the
observation {\bf 2} to be $n_1=n_2+k=n_3=n_4+k$.

This shows that $K_1$ and
$K_2$ can not be distinguished by the Vassiliev invariants of Legendrian
knots provided that they realize isotopic Pseudo-Legendrian knots. \qed
\end{emf}

\subsection{Some Proposition needed for 
the Proof of Theorem~\ref{equivalent}}

\begin{defin}\label{stabilization}

There are four types of small kinks that can be added to a Pseudo-Legendrian
knot. They are kinks of types $(1,1), (1,-1), (-1,-1), (-1,1)$, shown in
Figure~\ref{stabilization1.fig}. (Here the first number in the pair
corresponds to the increment to the rotation number of a knot and the second
to the increment to the selflinking number of the knot that occurs if one
adds a kink to a Pseudo-Legendrian knot in $\R^3$.)

\begin{figure}[htbp]
 \begin{center}
  \epsfxsize 10cm
  \hepsffile{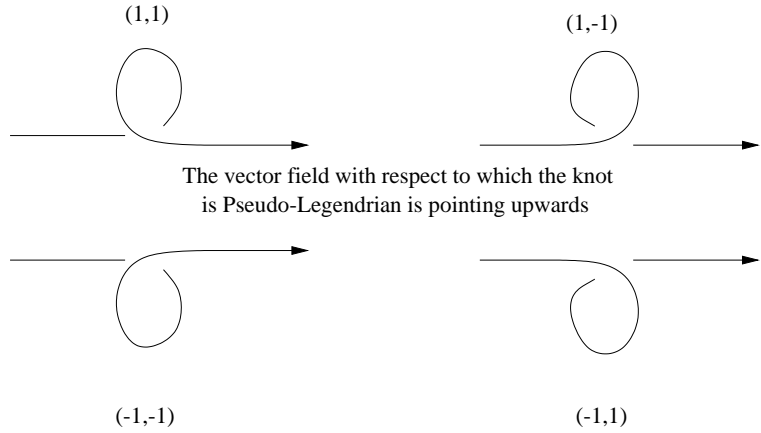}
 \end{center}
\caption{The four types of kinks}\label{stabilization1.fig}
\end{figure}

Let $K_p$ be a Pseudo-Legendrian knot. For $i,j\in\Z$ 
the $(i,j)$-{\em stabilization\/} $K_p^{i,j}$ of $K_p$ is a
Pseudo-Legendrian knot obtained from $K_p$ by the addition of 
$$\begin{cases}
i \text{ kinks of type } (1,1), & \text{ provided that } i\geq 0,\\
|i| \text{ kinks of type } (-1,-1), & \text{ provided that }i<0\\
\end{cases}$$
and of 
$$\begin{cases}
j  \text{ kinks of type } (-1,1), & \text{ provided that } j\geq 0,\\
|j| \text{ kinks of type } (1,-1), & \text{ provided that }j<0.\\ 
\end{cases}$$

Similarly one defines the $(i,j)$-{\em stabilization of a singular
Pseudo-Legendrian knot\/} with $n$ transverse double points.

It is easy to verify that the Pseudo-Legendrian 
isotopy type of the knot $K_p^{i,j}$ does not depend on the places
on $K$ where the kinks are added. (To get the same fact for singular
Pseudo-Legendrian knots one observes that it is possible to pull a kink
through a transverse double point of a singular Pseudo-Legendrian knot by a
Pseudo-Legendrian isotopy of the singular knot.)

We say that a Pseudo-Legendrian knot $K_p$ {\em admits finitely many
symmetric stabilizations\/} if there exists nonzero $n\in\Z$ such that
$K_p^{n,n}$ and $K_p$ are Pseudo-Legendrian isotopic. 
Similarly we introduce the notion of singular Pseudo-Legendrian knot
$K_{ps}$ that {\em admits finitely many
symmetric stabilizations.\/}
(Proposition~\ref{finitestab} implies that, if $K_p$ admits
finitely many symmetric stabilizations then for every $i,j\in\Z$ the knot $K_p^{i,j}$
also admits finitely many symmetric stabilizations.)

\end{defin}

\begin{figure}[htbp]
 \begin{center}
  \epsfxsize 10cm
  \hepsffile{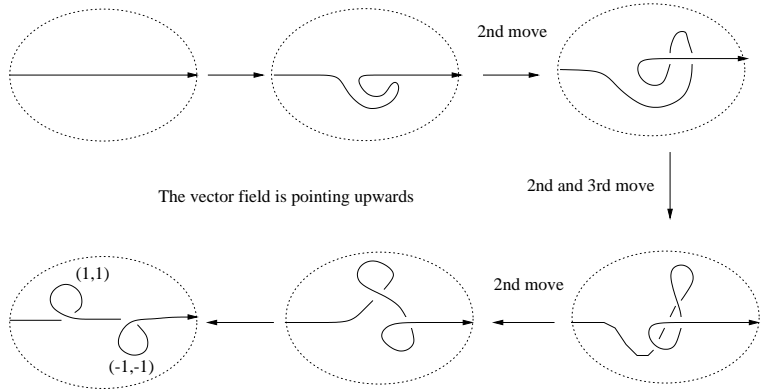}
 \end{center}
\caption{The creation of a $(1,1)$ and $(-1,-1)$ kinks by a
Pseudo-Legendrian isotopy}\label{kink2.fig}
\end{figure}

As it is shown in Figure~\ref{kink2.fig} the pair of kinks of types $(1,1)$
and $(-1,-1)$ can be cancelled by a Pseudo-Legendrian isotopy. Similar
considerations show that a pair of kinks of types $(1,-1)$ and $(-1,1)$ also
can be cancelled by a Pseudo-Legendrian isotopy.

Thus we get the following Proposition.

\begin{prop}\label{composestabilization}
For a Pseudo-Legendrian knot $K_p$ and for $i,j,k,l\in\Z$ the
$(i,j)$-stabilization $(K_p^{k,l})^{i,j}$ of $K_p^{k,l}$ is Pseudo-Legendrian
isotopic to $K_p^{i+k, j+l}$.
\end{prop}


\begin{prop}\label{stab}
Let $K_1$ and $K_2$ be Pseudo-Legendrian (with respect to $V_C$) knots that 
are $C^0$-close
isotopic to each other as unframed knots. Then there exist $i,j\in\Z$
such that $K_1$ is Pseudo-Legendrian isotopic to $K_2^{i,j}$. 
\end{prop}

\begin{emf}{\em Proof of Proposition~\ref{stab}\/}
Let $T$ be a tubular neighborhood of $K_1$ inside which $K_1$ and $K_2$ are
isotopic as unframed knots. 

Identify $(T,V_C\big|_T)$ with the solid torus in $\R^3$ such that both the 
axis of the torus and the vector field are parallel to the $z$-axis.
Similar to~\ref{description} we can depict the Pseudo-Legendrian knots
$K_1$ and $K_2$ by their knot diagrams obtained 
by projection to an annulus $A$ along the $z$-axis. 

Since $K_1$ and $K_2$ are isotopic inside $T$ as unframed knots, we get that
the knot diagram of $K_1$ in $A$ can be changed to a knot diagram of $K_2$
in $A$ by a sequence of Reidemeister moves. The second and the third
Reidemeister moves can be done in the Pseudo-Legendrian category. The first
Reidemeister move can not be done in the Pseudo-Legendrian category since
under it the velocity vector of a knot at one of the points becomes parallel
to the $z$-axis, and hence the knot becomes tangent to $V_C$ at one point.

There are four types of first Reidemeister move. They are distinguished by
the four possible types of kinks that appear under them, see
Figure~\ref{stabilization1.fig}. 

As it is shown in Figure~\ref{kink2.fig} one can create (or annihilate) a
$\{(1,1);(-1,-1)\}$-pair of kinks by a Pseudo-Legendrian isotopy. Similar
considerations show that one can create (or annihilate) a
$\{(1,-1);(-1,1)\}$-pair of kinks by a Pseudo-Legendrian isotopy.

Using this observation we can imitate the isotopy of unframed knots that
changes $K_1$ to $K_2$ by a Pseudo-Legendrian isotopy. 
Namely, if for example the isotopy that changes
$K_1$ to $K_2$ as an unframed knot involves the first Reidemeister 
that creates the $(1,1)$-kink, instead of this move we perform 
the Pseudo-Legendrian isotopy that creates a
$\{(1,1);(-1,-1)\}$-pair of kinks, make 
the $(-1,-1)$ kink very small and keep it very small during the 
later isotopy process.

In the end of this imitation process we get that $K_1$ is
Pseudo-Legendrian isotopic to a knot that looks exactly as $K_2$
except of a number of small extra kinks that are present on it. 

Cancel by a
Pseudo-Legendrian isotopy the $\{(1,1),(-1,-1)\}$ pairs of extra 
kinks either till there are no extra kinks of type $(1,1)$ left
or till there are no extra kinks of type $(-1,-1)$ left.

Cancel by a Pseudo-Legendrian isotopy 
the $\{(1,-1), (-1,1)\}$ pairs of extra kinks 
either till there are no
extra kinks of type $(1,-1)$ left or till there are no extra kinks of type 
$(-1,1)$ left.

It is clear that as the result of this process we obtain the
Pseudo-Legendrian isotopy of
$K_1$ to $K_2^{i,j},$ for some $i,j,\in\Z$.\qed
\end{emf}

\begin{prop}\label{finitestab}
Let $K_1$ and $K_2$ be Pseudo-Legendrian knots, and $i,j\in\Z$ be such that
$K_1$ and $K_2^{i,j}$ are Pseudo-Legendrian isotopic. Then for any
$k,l\in\Z$ the knots $K_1^{k,l}$ and $K_2^{i+k,j+l}$ are also
Pseudo-Legendrian isotopic.
\end{prop}
\begin{emf}{\em Proof of Proposition~\ref{finitestab}.\/}
Make the $|k|$ and $|l|$ extra kinks used to obtain $K_1^{k,l}$ very small
and concentrate them on a small piece of $K_1$. Keep them
small and close together during the isotopy process that was connecting
$K_1$ and $K_2^{i,j}$. As a result we get a Pseudo-Legendrian isotopy between 
$K_1^{k,l}$ and $(K_2^{i,j})^{k,l}$. Finally Proposition~\ref{composestabilization} 
says that $(K_2^{i,j})^{k,l}$ is Pseudo-Legendrian isotopic to
$K_2^{i+k,j+l}$. \qed.
\end{emf}

\begin{prop}\label{Legendrianstabilization}
Let $K$ be a Pseudo-Legendrian knot, that is Pseudo-Legendrian isotopic to a
Legendrian knot (i.e. its Pseudo-Legendrian 
isotopy type is realizable by a Legendrian knot). Then for any $i,j\in\N$ 
the Pseudo-Legendrian isotopy class of $K^{-i,-j}$ is also realizable by a
Legendrian knot.
\end{prop}
\begin{emf}{\em Proof of Proposition~\ref{Legendrianstabilization}.\/}
One observes that if $K_l$ is a Legendrian knot that is Pseudo-Legendrian
isotopic to $K$, then the Legendrian 
knot obtained by the modification shown in
Figure~\ref{twocusp.fig}~b is Pseudo-Legendrian isotopic to $K^{-1,0}$.
Similarly the Legendrian knot obtained by the modification shown in
Figure~\ref{twocusp.fig}~c is Pseudo-Legendrian isotopic to $K^{0,-1}$.
Performing the two modifications $i$ and $j$ times respectively we get the
Legendrian knot that is Pseudo-Legendrian isotopic to $K^{-i,-j}$.\qed
\end{emf}

\begin{prop}\label{atmostonecomponent}
Every component of the space of Pseudo-Legendrian curves contains at most 
one component of the space of Legendrian curves.
\end{prop}

{\em In Fact every component of the space of Pseudo-Legendrian curves contains
exactly one component of the space of Legendrian curves, but to prove this
statement, see Proposition~\ref{components} we have to use
Proposition~\ref{decrease} that is in turn based on
Proposition~\ref{atmostonecomponent}.\/}

\begin{emf}{\em Proof of Proposition~\ref{atmostonecomponent}.\/}
The $h$-principle says that the connected components of the space of
Legendrian curves in $(M,C)$ are in one to one correspondence with the
conjugacy classes of elements of $\pi_1(CM)$. Under this identification the
connected component of the space of Legendrian curves that contains a 
Legendrian curve $K_l$ corresponds to the conjugacy class of $\vec
K_l\in\pi_1(CM, \vec K_l(1))$, where the loop $\vec K_l$ is obtained by
mapping $t\in S^1$ to the point of $CM$ that corresponds to the velocity
vector of $K_l$ at $K_l(t)$.

Let $K_{0,l}$ and $K_{1,l}$ be Legendrian curves that are in $\mathcal L_p$. 
Consider a Pseudo-Legendrian homotopy $H:[0,1]\times S^1\rightarrow
M$ that connects $K_{1,l}$ and $K_{2,l}$. Observe that a Pseudo-Legendrian
curve $K$ also defines a loop $\vec K$ in $CM$ by mapping $t$ in $S^1$ to
the point of $CM$ that corresponds to the projection of the velocity vector
of $K$ at $K(t)$ to a contact plane along $V_C$. (Since $K$ is
Pseudo-Legendrian the projection is nonzero.)

Thus the family of Pseudo-Legendrian curves $K_t=H\big|_{t\times S^1}$
defines the free homotopy of loops $\vec K_{1,l}$ and $\vec K_{2,l}$. 
Hence by the $h$-principle $K_{1,l}$ and $K_{2,l}$ belong to the same component of the
space of Legendrian curves.

Since $K_{1,l}$ and $K_{2,l}$ were arbitrary Legendrian curves from
$\mathcal L$ we get that $\mathcal L_p$ contains at most one component of
the space of Legendrian curves. \qed
\end{emf}

\begin{prop}\label{everythingmadesingular} 
The statements of Propositions~\ref{composestabilization},~\ref{stab},
\ref{finitestab},~\ref{Legendrianstabilization} are true if one substitutes
everywhere in the statements of these Propositions Pseudo-Legendrian knots by
singular Pseudo-Legendrian knots, and Legendrian knots by singular
Legendrian knots.
\end{prop} 

The Proof of this Proposition is straightforward. One has to observe that it
is possible to pull a kink through a transverse double point of a singular
Pseudo-Legendrian knot by a Pseudo-Legendrian isotopy of the singular knot.

\begin{lem}\label{decrease}
Let $(M,C)$ be a contact manifold with a cooriented contact structure and
let $V_C$ be a vector field that coorients the contact structure.
Let $\mathcal L$ be a connected component of the space of
Legendrian curves in $(M,C)$ and let $\mathcal L_p$ be the corresponding
component of the space of Pseudo-Legendrian (with respect to $V_C$) curves. 
\begin{description}
\item[a]
Let $K\in\mathcal L_p$ be a Pseudo-Legendrian knot, then there exists 
$n\in \N$ such that the knot $K^{-n,-n}$ is Pseudo-Legendrian
isotopic to a Legendrian knot from $\mathcal L$.


\item [b] Let $K_s\in\mathcal L_p$ be a singular 
Pseudo-Legendrian knot (whose only singularities are $i$ transverse double
points), then there exists $n\in \N$ 
such that the knot $K_s^{-n,-n}$ is Pseudo-Legendrian 
isotopic to a singular Legendrian knot in $\mathcal L$.


\end{description}

\end{lem}

\begin{emf}{\em Proof of Lemma~\ref{decrease}.\/}
We prove statement {\bf a} of Lemma~\ref{decrease}. (The Proof of
statement {\bf b} of Lemma~\ref{decrease} is obtained in the same way as the
Proof of  statement {\bf a} of Lemma~\ref{decrease}, see
Proposition~\ref{everythingmadesingular}.)

The Theorem of Chow and Rashevskii~\cite{Chow},~\cite{Rashevskii} says that
that there exists a Legendrian knot $K_l$ that is $C^0$-small isotopic to
$K$ as an unframed knot. By Proposition~\ref{stab} we get that there exist
$i,j\in\Z$ such that $K^{i,j}$ is Pseudo-Legendrian isotopic to $K_l$.

Proposition~\ref{Legendrianstabilization} says that the Pseudo-Legendrian
isotopy classes of $K_l^{-1,0}$ and of $K_l^{0,-1}$ are also realizable by
Legendrian knots. By Propositions~\ref{composestabilization}
and~\ref{finitestab} $K^{i-1,j}$ is Pseudo-Legendrian isotopic to
$K_l^{-1,0}$ and $K^{i,j-1}$ is Pseudo-Legendrian isotopic to $K_l^{0,-1}$.
Using this fact we get that there exist $m\in\N$ such that $K^{-m,-m}$ is
Pseudo-Legendrian isotopic to a Legendrian knot.

The Pseudo-Legendrian homotopies shown in Figure~\ref{kink2.fig}
and~\ref{twokinksingular.fig} imply that $K^{-m,-m}$ is in $\mathcal L_p$.
Since by Propositions~\ref{atmostonecomponent} there is only one component of
the space of Legendrian curves contained in $\mathcal L_p$, we get that
$K^{-m,-m}$ is realizable by a Legendrian knot from $\mathcal L$. This
finishes the proof of statement {\bf a} of Lemma~\ref{decrease}.\qed

\begin{figure}[htbp]
 \begin{center}
  \epsfxsize 10cm
  \hepsffile{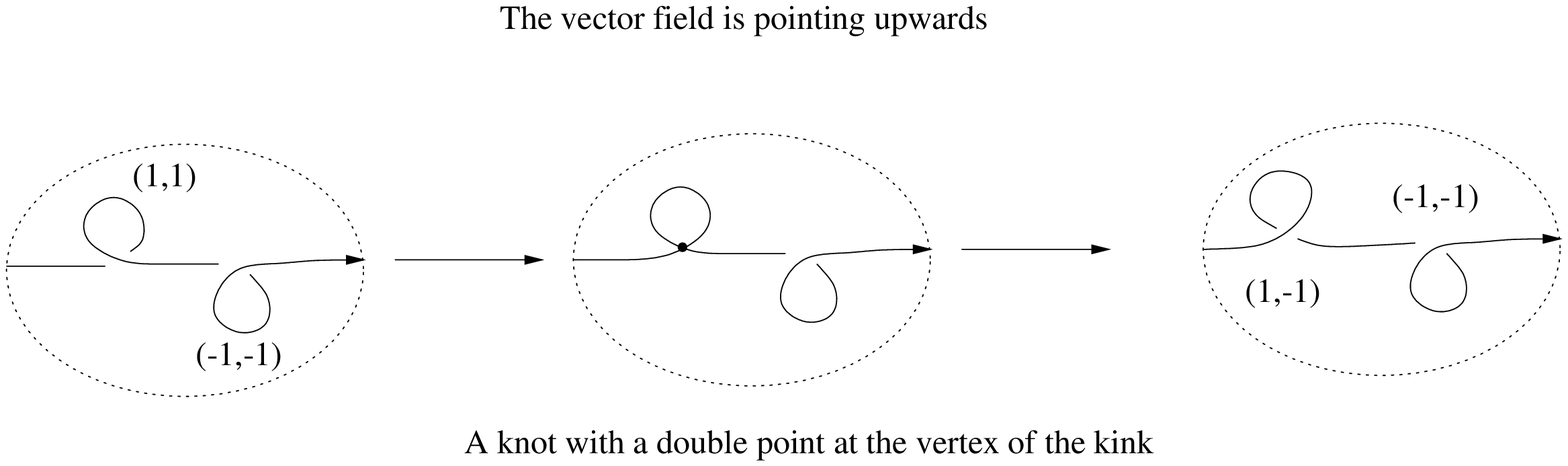}
 \end{center}
\caption{The creation of a $(1,1)$ and $(-1,-1)$ kinks by a
Pseudo-Legendrian isotopy}\label{twokinksingular.fig}
\end{figure}

\end{emf}

\begin{prop}\label{components}
Every component of the space of Pseudo-Legendrian curves contains exactly
one component of the space of Legendrian curves.
\end{prop}
\begin{emf}{\em Proof of Proposition~\ref{components}.\/}
Let $\mathcal L_p$ be a connected component of the space of
Pseudo-Legendrian curves, and let $K\in\mathcal L_p$ be a Pseudo-Legendrian
knot. In the proof of  statement {\bf a} of Lemma~\ref{decrease} 
we showed  that there exists
$n\in\N$ such that $K^{-n,-n}$ is realizable by a Legendrian knot. 
On the other hand the Pseudo-Legendrian isotopy shown in Figure~\ref{kink2.fig} 
and the Pseudo-Legendrian homotopy shown in Figure~\ref{twokinksingular.fig}
imply that $K^{-n,-n}$ is in $\mathcal L_p$. Thus $\mathcal L_p$ contains at
least one component of the space of Legendrian curves.
Proposition~\ref{atmostonecomponent} says that $\mathcal L_p$ contains at
most one component of the space of Legendrian curves. Thus it contains
exactly one component of the space of Legendrian curves.\qed
\end{emf}

\subsection{Proof of Theorem~\ref{equivalent}}\label{proofequivalent}
The fact that statement {\textrm I} of Theorem~\ref{equivalent} implies
statement {\textrm II} is clear.
Thus we have to  show that statement {\textrm II} implies statement
\textrm{\textrm I}. This is done by showing that there exists a homomorphism
$\psi:V_n^{\mathcal L}\rightarrow V_n^{\mathcal L_p}$ such that
$\phi\circ\psi=\id_{V_n^{\mathcal L}}$ and $\psi\circ\phi=\id_{V_n^{\mathcal
L_p}}$.

Let $x\in V_n^{\mathcal L}$ be an invariant.
In order to construct $\psi(x)\in V_n^{\mathcal L_p}$ we have to
specify the value of $\psi (x)$ on every Pseudo-Legendrian  $K\in\mathcal
L_p$.

\begin{emf}\label{definitionofpsi}{\em Definition of $\psi(x)$.\/} If the
isotopy class of the knot $K\in \mathcal
L_p$ is realizable by a Legendrian knot $K_l\in\mathcal L$, then put
$\psi(x)(K)=x(K_l)$. The value $\psi(x)(K)$ is well-defined because if
$K'_l\in\mathcal L$ is another knot realizing $K$, then
$x(K_l)=x(K'_l)$ by statement {\textrm I} of Theorem~\ref{isomorphism}.

Fix a Pseudo-Legendrian knot $K$. We explain how to define the value 
of $\psi(x)$ on all the $K^{n,n}$, $n\in\Z$.

If a Pseudo-Legendrian knot $K\in \mathcal L_p$ admits finitely many
symmetric stabilizations, then by Lemma~\ref{decrease} and Proposition 
~\ref{Legendrianstabilization} the knot 
$K$ and all the 
$K^{n,n}$ admit a Legendrian realization by a knot from $\mathcal L$. 
Thus we have defined the value of
$\psi(x)$ on all the Pseudo-Legendrian knots from $\mathcal L_p$ that 
admit finitely many symmetric stabilizations.

If $K\in \mathcal L_p$ admits infinitely many symmetric stabilizations, 
then either {\bf 1)} all the isotopy classes of $K^{q,q}$, $q\in\Z$, are
realizable by Legendrian knots from $\mathcal L$ or {\bf 2)} there exists
$q\in\Z$ such that $K^{q,q}$ 
is realizable by a Legendrian knot from $\mathcal L$ (see~\ref{decrease}) 
and $K^{q+1,q+1}$
is not realizable by a Legendrian knot from
$\mathcal L$. (In this case $K^{q+2,q+2}, K^{q+3,q+3}$ etc. 
also are not realizable by
Legendrian knots from $\mathcal L$, see~\ref{Legendrianstabilization}.) In the case {\bf
1)}
the value of $\psi(x)$ is already defined on
all the Pseudo-Legendrian knots from $\mathcal L_p$ that are
Pseudo-Legendrian isotopic to $K^{n,n}$, for some $n\in\Z$.

In the case {\bf 2)}
put
\begin{equation}\label{eqextension}
\psi(x)(K^{q+1,q+1})=\sum_{i=1}^{n+1}\bigl((-1)^{i+1}\frac{(n+1)!}{i!(n+1-i)!}
\psi(x)(K^{q+1-i, q+1-i})\bigr).
\end{equation}
(Proposition~\ref{decrease} implies that the sum on the right hand side
is well-defined.)
Similarly put
\[
\begin{array}{l}
 
\psi(x)(K^{q+2,q+2})=\sum_{i=1}^{n+1}\bigl((-1)^{i+1}\frac{(n+1)!}{i!(n+1-i)!}
\psi(x)(K^{q+2-i,q+2-i})\bigr), \\
\psi(x)(K^{q+3,q+3})=\sum_{i=1}^{n+1}\bigl((-1)^{i+1}\frac{(n+1)!}{i!(n+1-i)!}
\psi(x)(K^{q+3-i,q+3-i})\bigr)\text{ etc.}\\
\end{array}
\]

Now we have defined the value of $\psi(x)$ on all the Pseudo-Legendrian
knots from $\mathcal L_p$ that are symmetric stabilizations of a knot $K$ for
which case {\bf 2} holds.

Doing this for all the classes of $(n,n)$-stabilization equivalence for which
case {\bf 2} holds we define the value of $\psi(x)$
on all the knots from $\mathcal L_p$.
\end{emf}

{\bf Below we show that $\psi(x)$ is an order $\leq n$ invariant of framed
knots from
$\mathcal F$.}
We start by proving the following Proposition.
 
\begin{prop}\label{mainidentity}
Let $K^{q+1,q+1}$ be a Pseudo-Legendrian knot from $\mathcal L_p$, 
then $\psi(x)$ defined as
above satisfies identity~\eqref{eqextension}.
\end{prop}
 
\begin{emf}{\em Proof of Proposition~\ref{mainidentity}.\/}
If $K^{q+1,q+1}$ 
is not realizable by a Legendrian knot from $\mathcal L$, then
the statement of the proposition follows from the formula we used to
define $\psi(x)\bigl(K^{q+1,q+1}\bigr)$.
 
If $K^{q+1,q+1}$ is realizable by a Legendrian knot $K_l$, 
then consider a
singular Legendrian knot $K_{ls}$ with $(n+1)$ double points that are
vertices of $(n+1)$ small kinks such that we get $K_l$
if we resolve all the double points positively staying in the class of the
Legendrian
knots. (To create $K_{ls}$ we perform the first half of the
homotopy shown in Figure~\ref{kink.fig} in $n+1$ places on $K_l$.)

\begin{figure}[htbp]
 \begin{center}
  \epsfxsize 12cm
  \hepsffile{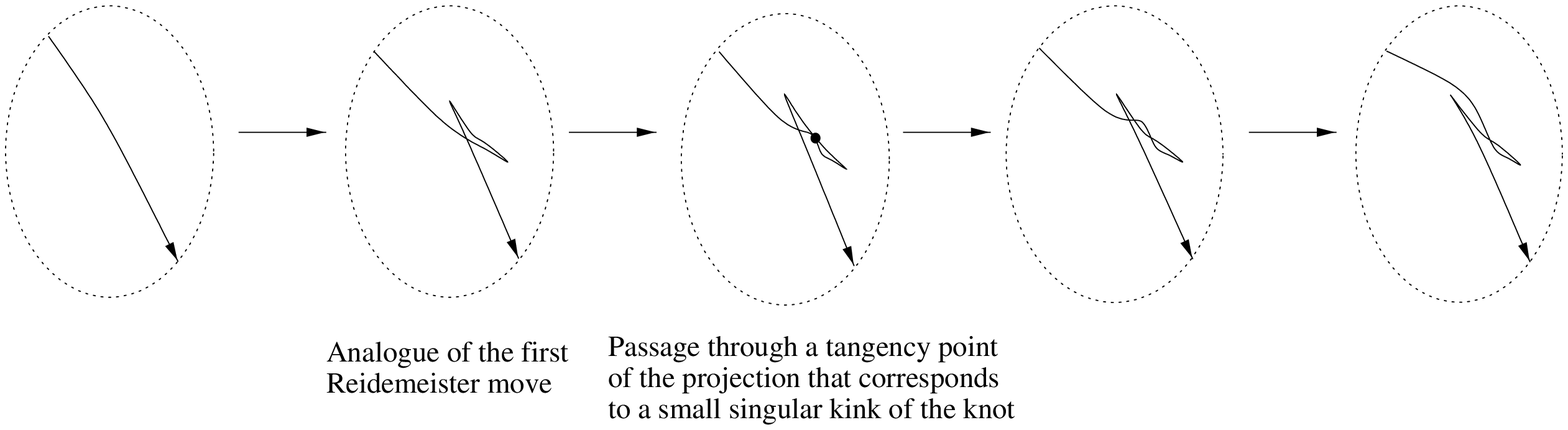}
 \end{center}
\caption{}\label{kink.fig}
\end{figure}

Let $\Sigma$ be the set of the $2^{n+1}$ possible resolutions of the double
points of $K_{ls}$.
For $\sigma\in \Sigma$ put $\sign(\sigma)$ to be the sign of the resolution,
and
put $K_{ls}^{\sigma}$ to be the nonsingular Legendrian knot obtained via the
resolution $\sigma$. Since $x$ is an order $\leq n$ invariant of Legendrian
knots we get that
\begin{multline}
0=\sum_{\sigma\in \Sigma}\bigl(\sign(\sigma)x(K_{ls}^{\sigma})\bigr)=\\
\psi(x)\bigl(K_l^{q+1,q+1}\bigr)+
\sum_{i=1}^{n+1}(-1)^i\frac{(n+1)!}{i!(n+1-i)!}\psi(x)\bigl
(K_l^{q+1-i,q+1-i}\bigr).
\end{multline}
(A straightforward verification shows that if 
we resolve $i$ double points of $K_{ls}$ negatively, then
we get the isotopy class of $K^{q+1-i,q+1-i}$.)
This finishes the proof of the Proposition. \qed
\end{emf}
 
\begin{emf}
Let $K_s\in\mathcal L_p$ be a singular Pseudo-Legendrian knot 
with $(n+1)$ double points.
Let $\Sigma$ be the set of the $2^{n+1}$ possible resolutions of the double
points of $K_s$. For $\sigma\in
\Sigma$ put $\sign(\sigma)$ to be the sign of the resolution, and put
$K^{\sigma}_s$
to be the isotopy class of the knot obtained via the resolution $\sigma$.
 
In order to prove that $\psi(x)$ is an order $\leq n$ invariant
of framed knots from $\mathcal L_p$, we have to show that
\begin{equation}\label{identitytoshow}
0=\sum_{\sigma\in \Sigma}\Bigl(\sign(\sigma)\psi(x)
\bigl(K_{s}^{\sigma}\bigr)\Bigr),
\end{equation}
for every singular $K_s\in \mathcal L_p$ with $(n+1)$ double points.
 
First we observe that the fact whether identity~\eqref{identitytoshow}
holds or not depends only on the Pseudo-Legendrian 
isotopy class of the singular Pseudo-Legendrian knot $K_s$
with $(n+1)$ double points.
 
If the isotopy class of $K_s$ is realizable by
a singular Legendrian knot from $\mathcal L$, then
identity~\eqref{identitytoshow} holds for $K_s$, since $x$ is an order $\leq
n$
invariant of Legendrian knots (and the value of $\psi(x)$ on a
Pseudo-Legendrian knot
$K\in\mathcal L_p$ realizable by a Legendrian knot $K_l\in\mathcal L$ 
was put to be $x(K_l)$).
 
Proposition~\ref{decrease} says that there exists $n\in\N$ such that 
the isotopy class of the singular Pseudo-Legendrian knot $K_s^{-n,-n}$
is realizable by a singular Legendrian knot from $\mathcal L$.

If $K_{s}$ admits finitely many symmetric stabilizations,
then all the isotopy classes of singular Pseudo-Legendrian knots 
from $\mathcal L_p$ that are isotopic to $K_s^{n,n}$, for some $n\in\Z$
are realizable by singular
Legendrian knots from $\mathcal L$, and we get that
identity~\eqref{identitytoshow} holds for $K_s$.
 
If $K_s$ admits infinitely
many symmetric stabilizations and all the isotopy classes of
singular Pseudo-Legendrian knots $K_s^{n,n}\in\mathcal L_p$, $n\in\Z$, 
are realizable by singular Legendrian knots from $\mathcal L$,
then~\eqref{identitytoshow} automatically holds for $K_s$.

If $K_s$ admits infinitely many symmetric stabilizations 
but not all the isotopy classes of $K_s^{n,n}$ $n\in\Z$, are
realizable by singular Legendrian knots from $\mathcal L$. 
Then let $q\in\Z$ be such that 
$K_s^{q,q}$ is realizable by a singular Legendrian knot from $\mathcal L$
and such that $K_s^{q+1,q+1}$ is not realizable by a singular Legendrian
knot from $\mathcal L$. (The existence of such $q$ follows from
Lemma~\ref{decrease}. Note that as it follows from the version of
Proposition~\ref{Legendrianstabilization} for singular knots,
see~\ref{everythingmadesingular}, the
Pseudo-Legendrian isotopy classes of $K_s^{q+2.q+2}, K_s^{q+3,q+3},
K_s^{q+4,q+4}$ etc. are also not realizable by singular Legendrian knot from
$\mathcal L$.)

Proposition~\ref{Legendrianstabilization} says that $K_s^{q-i,q-i}$, $i>0$, are realizable
by
singular Legendrian knots from $\mathcal L$ and hence
identity~\eqref{identitytoshow} holds for $K_s^{q-i,q-i}$, $i\geq 0$.
Using Proposition~\ref{mainidentity} and
the fact that identity~\eqref{identitytoshow} holds for
$K_s^{q-i,q-i}$, $i\geq 0$, we
show that~\eqref{identitytoshow} holds for $K_s^{q+1,q+1}$. Namely,
\begin{multline}
\sum_{\sigma\in
\Sigma}\sign(\sigma)\psi(x)\bigl({K_s^{q+1,q+1}}^{\sigma}\bigr)\\=
\sum_{\sigma\in
\Sigma}\Bigl(\sign(\sigma)
\sum_{i=1}^{n+1}(-1)^{i+1}\frac{(n+1)!}{i!(n+1-i)!}
\psi(x)\bigl((K_s^{(q+1-i),(q+1-i)})^{\sigma}\bigr)\Bigr)\\
=
\sum_{i=1}^{n+1}\Bigl((-1)^{i+1}\frac{(n+1)!}{i!(n+1-i)!}
\times\Bigl(\sum_{\sigma\in
\Sigma}\sign(\sigma)\psi(x)\bigl(({K_s^{(q+1-i), (q+1-i)}})^{\sigma}
\bigr)\Bigr)\Bigr)\\
=
\sum_{i=1}^{n+1}\Bigl(
(-1)^{i+1}\frac{(n+1)!}{i!(n+1-i)!}\times\Bigl(0\Bigr)\Bigr)=0.
\end{multline}

Similarly we show that ~\eqref{identitytoshow} holds for $K_{us}^{q+2,q+2},
K_{us}^{q+3,q+3}, \text{ etc}\dots$.
\end{emf}

\begin{emf} Clearly $\psi$ is a homomorphism and $\phi\circ\psi=\id_{V_n^{\mathcal L}}$.
 
Considering the values of $y\in V_n^{\mathcal L_p}$ on the $2^{n+1}$ possible
resolutions of a singular Pseudo-Legendrian knot with $n+1$ 
singular fragments shown in the middle part of Figure~\ref{twokinksingular.fig} we get
that $y$ should satisfy identity~\eqref{eqextension}.
Hence $\psi\circ\phi=\id_{V_n^{\mathcal
L_p}}$ and this finishes the proof of Theorem~\ref{isomorphism}. \qed
\end{emf}

{\bf Acknowledgments.} 
I am very grateful to Riccardo Benedetti, Yasha Eliashberg, Alexnadre
Kabanov, Stefan Nemirovskii, Christian Okonek, 
Sergei Tabachnikov, and Carlo Petronio for the their
motivating questions and suggestions as to what information about the
Legendrian knots could be captured using Vassiliev invariants.
And I am especially thankful 
to Riccardo Benedetti and Carlo Petronio who were, to my
knowledge, the first to conjecture that
a fact similar to the one we prove here is true.

These results were obtained during my stay at the Max-Planck-Institute f\"ur
Mathematik (MPIM), Bonn, 
and at the Institute for Mathematics, Z\"urich
University, and  
I would like to thank the Directors and the staff
of the MPIM, and the staff of the Z\"urich
University 
for hospitality and for providing the excellent working conditions.

\end{document}